# Initial version of State Transition Algorithm

Zhou Xiaojun, Yang Chunhua, Gui Weihua

Central South University, School of Information Science & Engineering, Changsha, Hunan, 410083, China

tiezhongyu2005@126.com

***Abstract*—In terms of the concepts of state and state transition, a new algorithm-State Transition Algorithm (STA) is proposed in order to probe into classical and intelligent optimization algorithms. On the basis of state and state transition, it becomes much simpler and easier to understand. As for continuous function optimization problems, three special operators named rotation, translation and expansion are presented. While for discrete function optimization problems, an operator called general elementary transformation is introduced. Finally, with 4 common benchmark continuous functions and a discrete problem used to test the performance of STA, the experiment shows that STA is a promising algorithm due to its good search capability.**

## I. INTRODUCTION

The concept of state means to a situation in which a material system remains, and it is featured in a group of quantities. The process of a system turning from a state into another one is called state transition, which can be described as a state transition matrix. The idea of state transition was created by a Russian mathematician named Markov when he expected to demonstrate a specific stochastic process (known as Markov process) [1]. Not only in communication theory but also in modern control theory, state transition matrix is of great importance. For example, in modern control theory, it can determine the stability of a system.

Based on the genetic and evolution mechanism of biology, the foraging behavior of ant colony, the rationale of metal annealing and the searching capability of bird flock, a variety of intelligent optimization algorithms have been proposed, which include genetic algorithm, ant colony algorithm, simulated annealing algorithm, particle swarm optimization algorithm and the rest, and they have made great contributions to different optimization problems. At present, there are numerous literatures researching in intelligent optimization algorithms, and they focus on the improvement on original algorithms and their hybrid applications as well as new algorithms' presenting [2].

The paper is to construct a new optimization method. Due to its foundation on state and state transition, the method is called State Transition Algorithm (STA), which has a simple form but has the ability to achieve various optimization problems. In the meanwhile, it is an open algorithm, in which any outstanding strategies can be incorporated, so as to reflect its strong capacity.

## II. UNDERSTANDING OF OPTIMIZATION ALGORITHMS

Considering the following unconstrained optimization problem

$$\min_{x \in R^n} f(x)$$

On the one hand, it is to analyze the problem in a classical view point [3,4].

In general, it adopts the iterative method. Let define $x_k$, $d_k$ and $a_k$ as the $k$th iteration point, direction of search and step, respectively, then the $k$th iteration can be described as

$$x_{k+1} = x_k + a_k d_k$$

The common way of selecting a step is by one dimensional search. While the methods of search direction include steepest descent method, conjugate gradient method, Newton method, quasi-Newton method, univariate search technique, Rosenbrock's rotating direction method, Powell method, and so on.

As described above, the classical method aims to search for a direction and a step. If classical algorithms are concerned in a state and state transition way, then an iterative point $x_k$ can be regarded as a state, the process of searching for a direction and a step will equate to a state transition process, and through a state transition, a new iterative point will be created.

On the other hand, it can also understand intelligent optimization algorithms in a state and state transition way. As for genetic algorithm[5,6], each individual can be considered as a state, and the updating population process of using genetic operators such as selection, crossover and mutation can equate to a state transition process. In the same way, for ant colony optimization[7,8], the ants adjusting their route by apperceiving pheromone, for simulated annealing[9,10], the solid determining its next state through the difference of internal energy and by the Metropolis criterion, and for particle swarm optimization[11-13], the bird flock updating its velocity and position can all be regarded as state transition processes.

## III. STATE TRANSITION ALGORITHM

With regard to the concept of state and state transition, a solution to the specific optimization problem can be described as a state, the thought of optimization algorithms (classical and intelligent optimization algorithms) can be treated as state transition, and the process to solve the optimization problem will become a state transition process.

Through the above analysis and discussion, it defines the following form of state transition

$$\begin{cases} x(k+1) = A_k x(k) + B_k u(k) \\ y_k = f(x(k+1)) \end{cases}$$

where, $x(k)$ stands for a state, corresponding to a solution to the optimization problem; then, $A_k$ and $B_k$ are state transition matrixes, which can be regarded as operators of optimization algorithm's thought; $u_k$ is the function with $x(k)$ and history states; and $f$ is the objective function.

### *A. Continuous function optimization problem*

Using various concepts of state space transformation for reference, it defines three special operators to solve the continuous function optimization problem.

(1) Rotation transformation

$$x(k+1) = (I_n + \alpha \times \frac{1}{n \times \|x(k)\|_2} \times R_r) \times x(k)$$

where, $x(k) \in \mathrm{R}^{n\times 1}$, $\alpha$ is a positive constant, called rotation factor. $R_r \in \mathrm{R}^{n\times n}$ , is a random matrix with its elements belonging to the range of [-1, 1] and$\| \cdot \|_2$ is 2-norm of a vector. Then, the following part will prove that the form has the function of rotation.

*Proof*:

$$\|x(k+1) - x(k)\|_2 = \left\| \alpha \times \frac{1}{n \times \|x(k)\|_2} \times R_r \times x(k) \right\|_2$$

$$= \frac{\alpha}{n\,\|x(k)\|_2} \|R_r \times x(k)\|_2$$

Taking the compatibility of matrix norm and vector norm into consideration, the $m_\infty$ norm in $\mathrm{R}^{n\times n}$ is compatible with the 2-norm in $\mathrm{R}^n$. Then

$$\|R_r \times x(k)\|_2 \le \|R_r\|_{m_\infty} \|x(k)\|_2$$

As a result,

$$\|x(k+1) - x(k)\|_2 \le \frac{\alpha}{n\,\|x(k)\|_2} \times \|R_r\|_{m_\infty} \|x(k)\|_2$$

$$= \frac{\alpha}{n} \times \|R_r\|_{m_\infty}$$

In addition, the elements of $R_r$ belong to the range of [-1, 1]; therefore, $\|R_r\|_{m_\infty} = n$, then

$$\|x(k+1) - x(k)\|_2 \le \alpha$$

It is easy to understand that the rotation definitely has the function of rotating. It aims to search in its vicinity, to be exactly, it explores solutions in a hyper sphere.

(2) Translation transformation

$$x(k+1) = x(k) + \beta \times R_T \times \frac{[x(k) - x(k-1)]}{\|x(k) - x(k-1)\|_2}$$

in which, $\beta$ is a positive constant, called translation factor. $R_T \in \mathrm{R}^1$ is a random variable, and its elements belong to the range of [0, 1].

We can understand that the function of translation transformation is to search along a line. As a matter of fact, it is also a way of one-dimensional search, to avoid some complicated computation.

(3) Expansion transformation

$$x(k+1) = x(k) + \gamma \times R_e \times x(k)$$

where, $\gamma$ is a positive constant, called expansion factor. $R_e \in \mathrm{R}^{n\times n}$ is a random diagonal matrix with its elements obeying the Gaussian distribution. It is well recognized that $y = \lambda x$ with $\lambda \neq 0$ can make a point $x$ of one-dimension expand in a positive or negative direction. To make a vector expand as freely as it can, the parameter $\lambda$ should become a diagonal matrix. It is easy to see that the expansion transformation has the function of searching in the whole space.

### *B. Discrete function optimization problem*

Different from the continuous part, under the discrete situation, a state also represents a solution, but it corresponds to a sequence. In a similar way, an important operator will be introduced.

In the subject of linear algebra, there is an important transformation matrix called elementary matrix, and the corresponding operation is called elementary transformation, which can swap two rows or two columns. The form of elementary matrix is in the following style

$$E(i,j) = \begin{bmatrix} 1 & & & & & & \\ & \ddots & & & & & \\ & & 0 & \cdots & 1 & & \\ & & \vdots & \ddots & & & \\ & & 1 & \cdots & 0 & & \\ & & & & & \ddots & \\ & & & & & & 1 \end{bmatrix} \begin{matrix} \\ \\ ith\ \ row \\ \\ jth\ \ row \\ \\ \\ \end{matrix}$$

It is manifest to find that the matrix originates from a unit matrix. As a matter of fact, some measures can be taken to make the style more useful. For example, we can exchange $m$ rows randomly to create a general elementary matrix $R$ on the base of unit matrix. The general elementary matrix is just the state transition matrix we called previously.

However, this is just the basic operator for discrete problem, how many rows should be exchanged and when should change the rows are really a key point deserve considering, that is to say, some other strategies are supposed to be introduced (not involved in the paper).

From the definition of state transition matrices for continuous and discrete optimization problems, the core thought of STA is introduced. Namely, different from other intelligent algorithms, a new solution (state) in STA is created by state transition transformation.

In STA, three main procedures are included: (1) Initialize parameters, create an initial state set. (2) Select the best state to implement state transitions, so as to create new states. (3) If the specified terminal criteria are met, it will stop; otherwise, it will continue to step (2).

## IV. EXPERIMENT AND SIMULATION

In order to evaluate the performance of STA, 4 common benchmark continuous functions and a traveling salesman problem (a classical discrete optimization problem) of 16 cities are introduced, as well as to illustrate the detailed procedures of the current version.

### *A. Test of continuous problems*

(1) Sphere function

$$f_1(x) = \sum_{i=1}^{n} x_i^2, -100 \le x_i \le 100$$

(2) Rosenbrock function

$$f_2 = \sum_{i=1}^{n} (100(x_{i+1} - x_i^2)^2 + (x_i - 1)^2), -30 \le x_i \le 30$$

(3) Rastrigin function

$$f_3 = \sum_{i=1}^{n} (x_i^2 - 10\cos(2\pi x_i) + 10), -5.12 \le x_i \le 5.12$$

(4) Griewank function

$$f_4 = \frac{1}{4000}\sum_{i=1}^{n} x_i^2 - \prod_{i=1}^{n} \cos\left|\frac{x_i}{\sqrt{i}}\right| + 1, -600 \le x_i \le 600$$

In current version of STA, the rotation factor will be decreased in an exponential way with base 4 from 1 to 1e-4 in an inner loop, that is to say, a rotation will only be halted when α is lower than 1e-4 in an iteration. And the translation factor and expansion factor are both fixed at 1.

Additionally, the times of implementing a transformation is called search enforcement. By the way, the translation will only be performed when a better state is gained, and all transformations will only accept the best of the state sets. In the experiment, the search enforcement is 32, the dimension size of all problems will be 10, 20 and 30, and the corresponding maximum iteration numbers are 1000, 1500 and 2000, respectively. There are 50 independent trails carried out, and the statistical results are given in Table 1.

Table 1 Test results of continuous functions

| Fcn | Dim | Iterations | Performance | | |
|---|---|---|---|---|---|
| | | | best | mean | std |
| $f_1$ | 10 | 1000 | 0 | 0 | 0 |
| | 20 | 1500 | 0 | 0 | 0 |
| | 30 | 2000 | 0 | 0 | 0 |
| $f_2$ | 10 | 1000 | 0.0607 | 2.8419 | 16.1679 |
| | 20 | 1500 | 0.1734 | 4.1075 | 1.0313 |
| | 30 | 2000 | 2.5937 | 16.1735 | 15.0863 |
| $f_3$ | 10 | 1000 | 0 | 9.5544 | 13.1377 |
| | 20 | 1500 | 0 | 60.5131 | 40.2447 |
| | 30 | 2000 | 38.8133 | 158.1179 | 46.4595 |
| $f_4$ | 10 | 1000 | 0 | 0.1512 | 0.1741 |
| | 20 | 1500 | 0 | 0.0293 | 0.0701 |
| | 30 | 2000 | 0 | 0.0025 | 0.0128 |

As shown in Table 1, the STA has good performance for $f_1$. For $f_2$, it is a little deficient because it cannot obtain the global optimum; however, the mean is acceptable in some degree. With regard to $f_3$, the STA has the ability to find the global optimum; nevertheless, the probability is not very high due to the increase of the mean. Concerning the decrease of the mean for $f_4$, it is satisfactory to say that the STA is advantageous for functions of this kind.

### *B. Test of a discrete problem*

Taking a general elementary matrix $R \in \mathbb{R}^{5\times 5}$ for example, its form includes all of the following styles but not restricted to these.

$$\begin{bmatrix} 1 & 0 & 0 & 0 & 0 \\ 0 & 1 & 0 & 0 & 0 \\ 0 & 0 & 0 & 0 & 1 \\ 0 & 0 & 0 & 1 & 0 \\ 0 & 0 & 1 & 0 & 0 \end{bmatrix}, \begin{bmatrix} 0 & 0 & 0 & 1 & 0 \\ 0 & 1 & 0 & 0 & 0 \\ 0 & 0 & 1 & 0 & 0 \\ 0 & 0 & 0 & 0 & 1 \\ 1 & 0 & 0 & 0 & 0 \end{bmatrix}, \begin{bmatrix} 0 & 0 & 0 & 1 & 0 \\ 0 & 1 & 0 & 0 & 0 \\ 1 & 0 & 0 & 0 & 0 \\ 0 & 0 & 0 & 0 & 1 \\ 0 & 0 & 1 & 0 & 0 \end{bmatrix}$$

If a initial sequence is [1,2,3,4,5], then after the above general elementary matrices, it will becomes [1,2,5,4,3], [4,2,3,5,1], [4,2,1,5,3] respectively. It is clearly to understand that after a general elementary transformation, and a sequence can change its order arbitrarily, which just meets the demand for global search and local search.

The successive part will illustrate the steps of solving a discrete function optimization problem.

Step1: Create some sequences randomly and uniformly as initial state set.

Step2: Select the best state from the initial state set. At this point, the best sequence is [ 10 11 6 12 7 13 14 9 16 8 5 2 3 4 1 15 10], with its path length at 113.6531 units and the order of traversing in Fig(a). Then, some appropriate general elementary transformations will be operated on the best sequence. After one operation, the order of traversing the cities is shown in Fig(b).

Step3: Repeat some procedures in step2 until the terminal criteria are met. The final order of traversing the cities is illustrated in Fig(c).

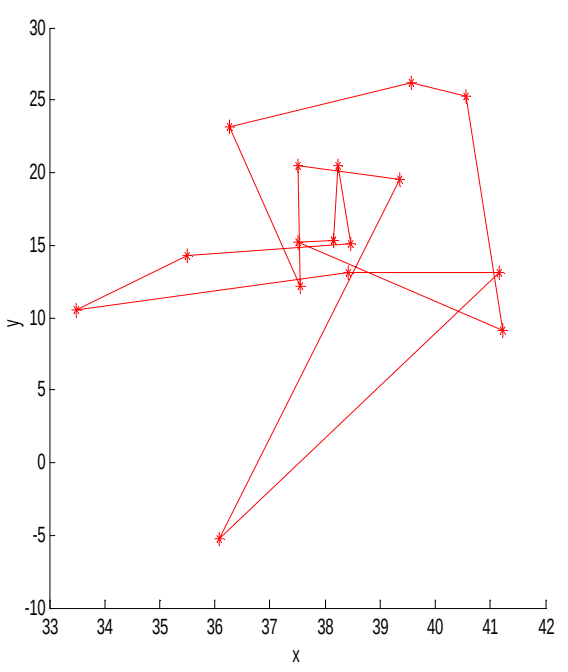

Fig(a). order in initiation

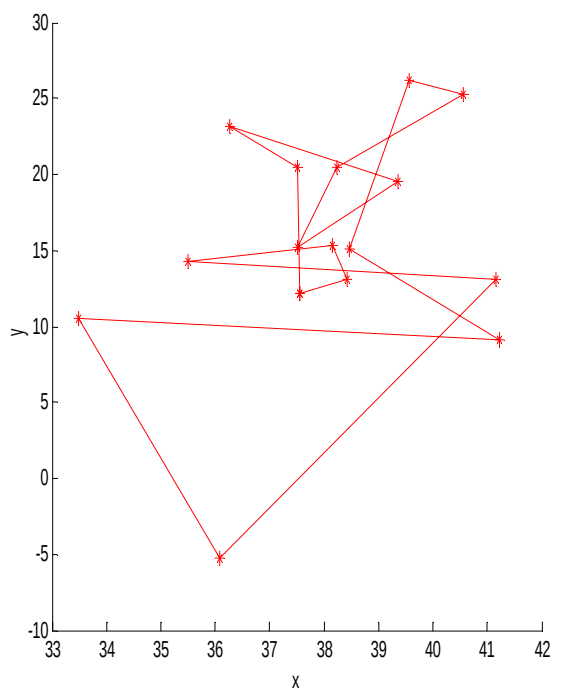

Fig(b). order after an operation

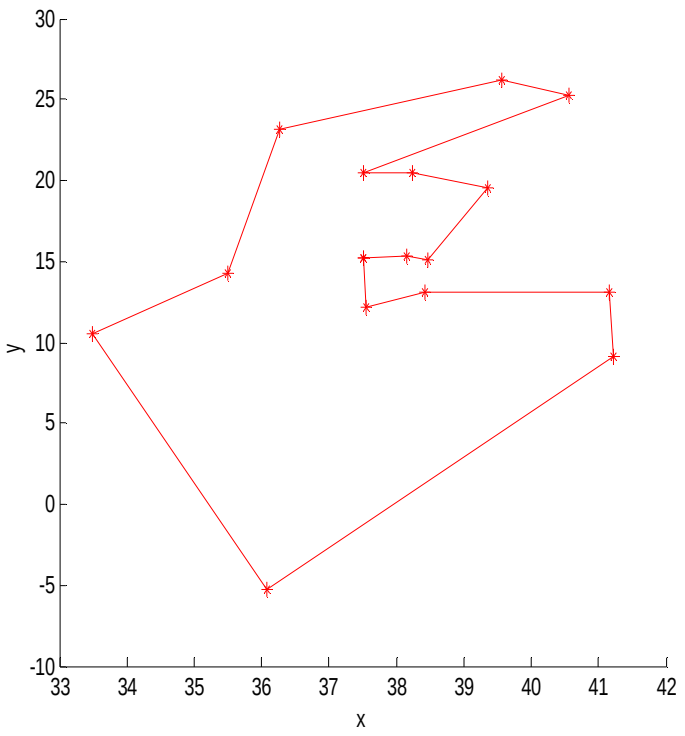

Fig(c). order under terminal criteria

## V. CONCLUSION

Based on the state and state transition, the STA, not only has the simple form but also possess clear physical connotation. For continuous function optimization problems, it presents the state transitions including rotation transformation, translation transformation and expansion transformation; as for discrete function optimization problems, it proposes the general elementary transformation matrix, which means that STA has strong adaptability to different problems. Using the 4 benchmark continuous functions and a traveling salesman problem for test, results show that STA has good search ability.

In addition, for one thing, the current version of STA, it utilizes the rotation transformation for local search, and the expansion transformation for global search, while the translation transformation is used just for simplifying the one-dimension search as well as balancing them. In comparison with other intelligent algorithms, although STA seems to be not as intelligent as them (most intelligent algorithms originate from the living principle of creatures in nature), STA has provided a new method for optimization. Also as a stochastic search method, to some extent, STA seems an algorithm not based on population, because the transformation is performed only on the best state. In respect of the mechanism of learning, STA only has self-learning function without any communication. For another, because the current version of STA is just a method of updating of individual, STA can possess other styles. For example, it can transform the states in a group way, in spite of their complex operations. In any case, it indicates that STA has much room for development, and it will be a promising method in optimization.

## ACKNOWLEDGMENT

The work was supported by the National Science Found for Distinguished Young Scholars of China (Grant No. 61025015).